\documentclass[10pt]{amsart}
\usepackage{amsmath}
\usepackage{amscd}
\usepackage{graphics}
\theoremstyle{plain}

\newtheorem{Thm}{Theorem}

\theoremstyle{definition}
\newtheorem{Rem}{Remark}

\def\IZ{\mathbb Z}

\def\cpb{\hbox{$\overline{{\mathbb C}P^2}$}}

\def\s2x{\hbox{$S^2 \times S^2$}}

\def\S{\Sigma}

\newenvironment{prooff}{\medskip \par \noindent {\it Proof}\ }{\hfill
$\square$ \medskip \par}
	\def\sqr#1#2{{\vcenter{\hrule height.#2pt
   		\hbox{\vrule width.#2pt height#1pt \kern#1pt
      		\vrule width.#2pt}\hrule height.#2pt}}}
	\def\square{\mathchoice\sqr67\sqr67\sqr{2.1}6\sqr{1.5}6}

\begin{document}
\title[]{Noncomplex smooth 4-manifolds with genus-2
Lefschetz fibrations} 
\author{Burak Ozbagci \and Andr\'as I. Stipsicz}
\address{Department of Mathematics\\University of California\\
Irvine, CA 92697 \and Department of Analysis, ELTE TTK, 
M\'uzeum krt. 6-8, Budapest, Hungary}
\email{bozbagci@math.uci.edu \and stipsicz@cs.elte.hu}
\date{\today}
\begin{abstract}
We construct noncomplex smooth 4-manifolds which admit 
genus-2 Lefschetz fibrations over $S^2$. The fibrations are
necessarily hyperelliptic, and the resulting 4-manifolds are not 
even homotopy equivalent to complex surfaces. Furthermore, these examples
show that fiber sums of holomorphic Lefschetz fibrations do not necessarily
admit complex structures.
\end{abstract}
\maketitle

\setcounter{section}{-1}
In the following we will prove the following theorem.

{\Thm
There are infinitely many (pairwise nonhomeomorphic) 4-manifolds
which admit genus-2 Lefschetz fibrations but do not
carry complex structure with either orientation.}

\bigskip

\begin{figure}[h]
 \begin{center}
    \includegraphics{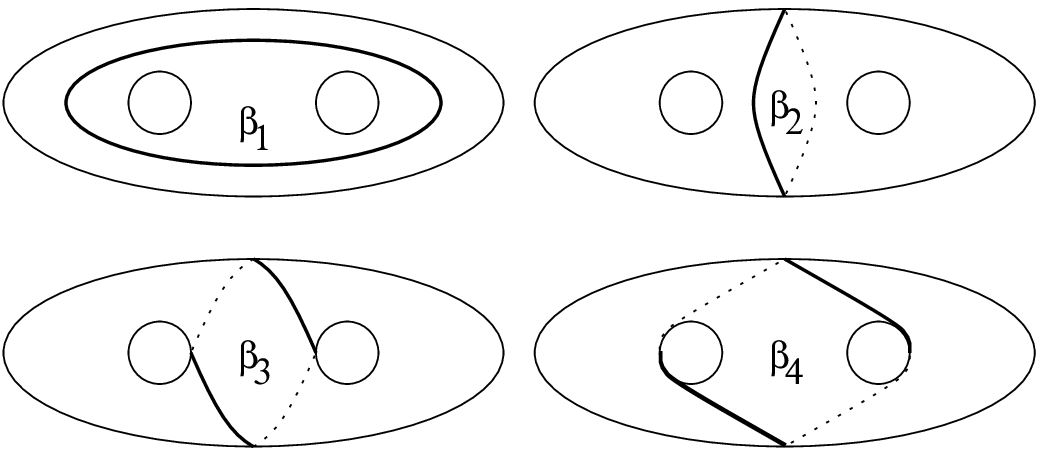}
  \caption{} \label{matsu}
   \end{center}
 \end{figure}

Matsumoto \cite{m}  showed that $S^2 \times T^2 \# 4 \cpb$ admits a 
genus-2 Lefschetz fibration over $S^2$ with global monodromy
$({\beta}_1 ,...,{\beta}_4 ) ^2$, where ${\beta}_1 ,...,{\beta}_4$ 
are the curves indicated by Figure~\ref{matsu}.
(For definitions and details regarding Lefschetz fibrations
see \cite{m}, \cite{gs}.)

Let $B_n$ denote the smooth 4-manifold which admits a 
genus-2 Lefschetz fibration over $S^2$ 
with global monodromy 
$$(({\beta}_1 ,...,{\beta}_4)^2 ,(h^n({\beta}_1),...,h^n({\beta}_4))^2) $$
where $h=D(a_2)$ is a positive Dehn twist about the curve $a_2$ 
indicated in Figure ~\ref{pi}. 

{\Thm
For the 4-manifold $B_n$ given above we have ${\pi}_1 (B_n)= 
{\IZ} \oplus {\IZ}_n $.}  

\begin{proof}

Standard theory of Lefschetz fibrations gives that
$${\pi}_1 (B_n) = {\pi}_1 ({\S}_2) / < {\beta}_1 ,...,{\beta}_4 , 
h^n ({\beta}_1),..., h^n ({\beta}_4)> .$$
Let $\{ a_1 ,b_1 ,a_2 ,b_2 \}$ be the standard generators 
for ${\pi}_1 ({\S}_2)$ (Figure ~\ref{pi}). 

\begin{figure}[h]
 \begin{center}
    \includegraphics{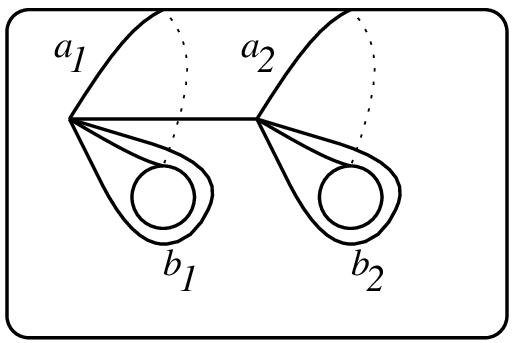}
  \caption{} \label{pi}
   \end{center}
 \end{figure}

Then we observe that 
 
${\beta}_1 = b_1 b_2$,

${\beta}_2 = a_1 b_1 a_1^{-1} b_1^{-1} = a_2 b_2 a_2^{-1} b_2^{-1}$,

${\beta}_3 = b_2 a_2 b_2^{-1} a_1$,

${\beta}_4 = b_2 a_2 a_1 b_1$,

$h^n ({\beta}_1) = b_1 b_2 a_2^n$,

$h^n ({\beta}_2) = {\beta}_2$,

$h^n ({\beta}_3) = {\beta}_3$,

$h^n ({\beta}_4) = b_2 a_2^{n+1}  a_1 b_1$.

Hence 

${\pi}_1 (B_n) = < a_1 , b_1 , a_2 , b_2 \ | \ b_1 b_2 , 
\ [a_1 , b_1] ,\ [a_2 , b_2] , \ b_2 a_2 b_2^{-1} a_1 ,
\ b_2 a_2 a_1 b_1 ,\  b_1 b_2 a_2^n ,\ b_2 a_2^{n+1}  a_1 b_1 >$ 

$=<a_2 , b_2 \; | \; [a_2 , b_2] , \ a_2^n >= \IZ \oplus {\IZ}_n,$
and this concludes the proof.
\end{proof}

{\Thm 
$B_n$ does not admit a complex structure.}

\begin{proof}

Assume that $B_n$ admits a complex structure. 
Let $M_n$  denote its 
$n$-fold cover  for which $\pi _1(M_n)\cong \IZ$
and $M_n'$ the minimal model of $M_n$. By the theorem
of Gompf \cite{gs} $B_n$ admits a symplectic structure, hence 
so does $M_n$ and (by combining results of Taubes and Gompf 
\cite{t}, \cite{g2}) $M_n'$. Consequently,  if $B_n$ is a complex
surface, then we have a symplectic, minimal complex surface $M_n'$ with 
$\pi _1(M_n')\cong \IZ$. In the following we will show that this leads to
a contradiction.

By the Enriques-Kodaira 
classification of complex surfaces \cite{b}, (since $b_1(M_n')=1$)
$M_n'$ is either
a surface of class $VII$ (in which case $b_2^+ (M_n')=0$),
a secondary  Kodaira surface (in which case $b_2(M_n')=0$)
or a (minimal) properly elliptic surface.

Since $M_n'$ is a symplectic 4-manifold, 
$b_2^+(M_n')$ (and so $b_2(M_n')$) is positive; this observation excludes
the first two possibilities. 

Suppose now
that $M_n'$ admits an elliptic fibration
over a Riemann surface.
If the Euler characteristic of $M_n'$ is 0, then (following form
the fact that $b_1(M_n')=b_3(M_n')=1$) we get that $b_2(M_n')=0$,
which leads to the above contradiction.
Suppose finally that $M_n'$ is a  minimal elliptic surface with positive
Euler characteristic.
Since $b_1(M_n')=1$, it can only be fibered over $S^2$ (see for example
\cite{fm}). 
In that case (according to \cite{g}, for example) its fundamental group is
$$\pi _1(M_n')= < x_1,\ldots x_k \  |\ x_i ^{p_i}=1, \ i=1,\ldots , k;\ 
x_1\cdot \cdot \cdot x_k =1>.$$
This cannot be isomorphic to $\IZ$, since if $\pi _1(M_n)\cong \IZ =<a>$,
then $x_1=a^{m_1}$ for some $m_1\in \IZ$, 
so $a$ has finite order, which is a contradiction.
Consequently the assumption that $B_n$ is complex leads us to a contradiction,
hence the theorem is proved.
\end{proof}

{\Rem
The above proof, in fact, shows that $B_n$ is not even 
homotopy equivalent to a complex surface --- our arguments used 
only homotopic invariants (the fundamental group, $b_2$ and $b_2^+$)
of the 4-manifold $B_n$. Note that basically the same idea shows
that ${\overline {B}}_n$ (the manifold $B_n$ with the opposite 
oreintation) carries no complex structure: The arguments involving
the fundamental group, $b_2$ and the Euler characteristic only, apply
without change. Since the fiber of the Lefschetz fibration on $B_n$
is homotopically essential and provides a class
with square 0, the intersection form of $B_n$ and 
so of $M_n$ are not definite --- consequently these manifolds cannot
be homotopy equivalent (with either orientation)
to the blow-up of a surface of Class $VII$.}

\begin{prooff}{\em of Theorem 1.}
By the definition of the 4-manifolds $B_n$ we get 
infinitely many  manifolds admitting genus-2 (consequenlty
hyperelliptic) Lefschetz fibrations which are (by Theorem 2.)
nonhomeomorphic. As Theorem 3. and the above remark show, the manifolds
$B_n$ do not carry complex structures with either orientation, hence
the proof of the Theorem 1. is complete.
\end{prooff}

{\Rem
We would like  to point out that similar examples have been found by 
Fintushel and Stern \cite{fs} --- they used Seiberg-Witten theory 
to prove that their (simply connected) genus-2
Lefschetz fibrations  are noncomplex.}

Note  that $B_n$ is given as the fiber sum of two copies of
$S^2\times T^2 \# 4\cpb$, hence provides an example of the phenomenon
that the fiber sum of holomorphic Lefschetz fibrations is not necessarily 
complex. 

{\it Acknowledgement.} \ Examples of genus-2 Lefschetz fibrations 
with $\pi_1 = \IZ \oplus {\IZ}_n$ were also 
constructed (as fiber sums) independently by Ivan Smith \cite{s}.


\begin{thebibliography}{99999}

\bibitem[BPV]{b}
W. Barth, C. Peters and A. Van de Ven, {\em Compact complex surfaces,} 
Springer-Verlag, 1984.

\bibitem[FM]{fm} R. Friedman and J. Morgan, 
{\em Smooth four-manifolds and complex surfaces,}
Springer-Verlag, 1994.

\bibitem[FS]{fs}
R. Fintushel and R. Stern, {\em Private communication.}

\bibitem[G1]{g} 
R. Gompf, {\em Nuclei of elliptic surfaces,} Topology {\bf30} (1991), 
479--511. 

\bibitem[G2]{g2}
R. Gompf, {\em A new construction of symplectic manifolds,}
Ann. of Math. {\bf142} (1995), 527--598.

\bibitem[GS]{gs} 
R. Gompf and A. Stipsicz,
{\em An introduction to 4-manifolds and Kirby calculus,}
book in preparation.

\bibitem[M]{m}
Y. Matsumoto, {\em Lefschetz fibrations of genus two - a topological
approach,} Proceedings of the 37th Taniguchi Symposium on
Topology and Teichm\"{u}ller Spaces, ed. Sadayoshi Kojima
et al., World Scientific (1996), 123--148.

\bibitem[S]{s} I. Smith, {\em Ph.D thesis.} 

\bibitem[T]{t} 
C. Taubes, {\em $SW\Longrightarrow Gr$: From the Seiberg-Witten
equations to pseudo-holomorphic curves,}
Journal of the Amer. Math. Soc. {\bf9} (1996), 845--918.

\end{thebibliography}
\end{document}